\documentclass[pdflatex,sn-mathphys]{sn-jnl}

\usepackage{graphicx}
\usepackage{float}
\usepackage{array}
\usepackage{tabularx}
\usepackage{amsmath}
\usepackage{enumerate}
\usepackage{xcolor}
\usepackage{soul} 
\usepackage{hyperref}
\hypersetup{
    colorlinks=true,
    linkcolor=blue,
    filecolor=magenta,      
    urlcolor=black,
    pdftitle={Overleaf Example},
    pdfpagemode=FullScreen,
    }

\theoremstyle{thmstyleone}%
\newtheorem{theorem}{Theorem}

\DeclareMathOperator{\stt}{s.t.} 

\DeclareMathOperator*{\argmin}{argmin}

\def\R {{\mathbb R}}

\graphicspath{ {./figures/} }

\begin{document}
\title{Cartoon-texture evolution for two-region image segmentation}
%
\author*[1]{\fnm{Laura} \sur{Antonelli}}\email{laura.antonelli@cnr.it}
\equalcont{These authors contributed equally to this work.}

\author[2]{\fnm{Valentina} \sur{De Simone}}\email{valentina.desimone@unicampania.it}
\equalcont{These authors contributed equally to this work.}

\author[2]{\fnm{Marco} \sur{Viola}}\email{marco.viola@unicampania.it}
\equalcont{These authors contributed equally to this work.}

\affil*[1]{\orgdiv{Institute for High Performance Computing and Networking (ICAR)}, \orgname{National Research Council (CNR)}, \orgaddress{\street{via Pietro Castellino, 111}, \city{Naples}, \postcode{80131}, \country{Italy}}}

\affil[2]{\orgdiv{Department of Mathematics and Physics}, \orgname{University of Campania ``Luigi Vanvitelli''}, \orgaddress{\street{viale Abramo Lincoln, 5}, \city{Caserta}, \postcode{81100}, \country{Italy}}}


%
\abstract{
Two-region image segmentation is the process of dividing an image into two regions of interest, i.e., the foreground and the background. 
To this aim, Chan et al. \cite{bib:CEN2006} designed a model well suited  for smooth images.
One drawback of this model is that it may produce a bad segmentation when the image contains oscillatory components.
Based on a cartoon-texture decomposition of the image to be segmented,  
we propose a new model that is able to
produce an accurate segmentation of images also containing noise or oscillatory information like texture. The novel model
leads to a non-smooth constrained optimization problem which we solve by means of the ADMM method. The convergence of the numerical scheme is also proved. Several experiments on smooth, noisy, and textural images show the effectiveness of the proposed model.}

\keywords{Image segmentation, cartoon-texture decomposition, non-smooth optimization, ADMM method}

\maketitle              
\section{Introduction}
Image segmentation  is a fundamental task in image processing and computer vision. 
It consists in dividing an image into non-overlapping regions of shared features, such as intensity, smoothness, and texture, which are related to the final goal of the segmentation. Thus, the division into regions is not unique, and the image segmentation can be regarded as a strongly ill-posed problem.\\
Let $f$ be an image defined in a domain $\Omega \subset \R^d$ ($d \geq 2$),
segmenting $f$ 
consists in finding a  decomposition of the domain {\small $\Omega$} into a set of non-empty pairwise-disjoint regions $\Omega_i$, $i=1,\ldots,m$.
\noindent A segmentation of $f$ can be expressed through a curve $C^*$ that matches the boundaries of the decomposition of {\small $\Omega$}, i.e. {\small$C^*= \bigcup \limits_{i =1 }^m\partial \Omega_i$} and/or   a piecewise-constant function $f^*$  defined on {\small $\Omega$} that approximates $f$.\\
The research on image segmentation has
made several advances in the last decades and various approaches have been developed, 
including thresholding, region
growing, edge detection and variational methods. \cite{bib:sur,bib:sur2,bib:sur3}.
Variational models, based on optimizing energy functionals, have been widely investigated, proving to be very effective
on different images;
curve evolution~\cite{bib:sethian}, anisotropic
diffusion~\cite{bib:PeronaMalik} and the Mumford-Shah model~\cite{bib:MumfordShah1989} are good representatives of these
methods. 
Other recent approaches to image segmentation include learning-based methods, which often exploit
deep-learning techniques~\cite{bib:dl0,bib:dl2,bib:dl3}, watershed \cite{bib:Challa2019}, random walk methods \cite{bib:Aletti2021}, graph cuts \cite{bib:Niazi2022,bib:He2019}, epidemiological models on images \cite{bib:Bampis2017}. However, in this case,
a large amount of data must be available to train learning networks, thus making those approaches impractical in
some applications.\\
Two-region segmentation is here considered, where the domain of the given image $ \bar{f}$ is separated in two regions of interest, so $m=2$ and $\Omega=\Omega_{in} \cup \Omega_{out}$, i.e. $\Omega_{in}$ and $\Omega_{out}$ are the foreground and the background of the image, respectively. Although the choice of $m=2$ significantly simplifies the segmentation problem, it has a lot of
application fields, such as biological and medical imaging, text
extraction, compression of screen content and mixed content documents, and can be used as a computational kernel for more complex segmentation tasks  \cite{bib:application1,bib:Zhang2008ExtractionOT,bib:7533087,bib:application2,bib:multiregions,bib:gregoretti2016}.

\noindent A widely-used two-region model was introduced by Chan and Vese in \cite{bib:ChanVese2001} and,  together with its variations, is regarded as state of the art in the segmentation community. These models are currently used in medical and astronomical application fields and have lately been associated with machine learning frameworks (see, e.g.~\cite{CEN1,CEN2,bib:dl1,bib:dl2,CEN5,CEN3}).
The Chan-Vese model is 
 a special case of the most popular Mumford-Shah one \cite{bib:MumfordShah1989} restricted to piecewise constant functions. The solution is the best
approximation to $\bar{f}$ among all the functions that take only two values, $c_{in}$ and $c_{out}$.   
As is the case of many variational models for image processing,
the model results in a non-convex optimization problem and may have various local minima. 
Chan, Esedo\=glu and Nikolova~\cite{bib:CEN2006} propose a convexed relaxation model, here denoted as CEN, which considers the case of $f$ taking values in $[0,1]$, and sets one of the two regions as
\begin{equation*}
    \Omega_{in} = \left\{ x : f(x) > \alpha \right\} \mbox{ for a.e. } \alpha \in (0,1).    
\end{equation*}
The CEN model first computes the values $c_{in}$ and $c_{out}$, and then, given $\lambda>0$, it determines $f$ by solving the convex minimization problem
\begin{equation} \label{eq:cen_continuous}
\begin{array}{rl}
\displaystyle \min_{0 \le
 f \le 1} & \displaystyle \int_{\Omega} \vert\nabla f\vert \, dx
                        + \lambda \int_{\Omega} \left( (c_{in} -  \bar{f} (x))^2 f(x) + (c_{out} -  \bar{f} (x))^2  (1-f(x)) \right) dx.
\end{array}
\end{equation}

\noindent We note that the aforementioned models assume that each image region is defined as a smooth or constant function. However, images may not be piecewise smooth or flat as a whole, but they may contain some non-smooth regions. In practice, imposing smoothness on such kind of images may lead to a destructive averaging of the image content \cite{bib:6566192}, which can produce an inaccurate segmentation. Exploiting information on the non-smooth structure of an image can help to improve the CEN model to be effective on a larger set than the one of smooth-images as done, e.g. in \cite{bib:Antonelli2020Adaptive}, thanks to the introduction of spatially-varying regularization methods.

In this paper we will design a new model for two-region image segmentation that, starting from a rough cartoon-texture decomposition $\bar{f}=\bar{u}+\bar{v}$ of the initial image, produces a cartoon-texture-driven decomposition of $\bar{u}$ as $\bar{u} = u + v$ and simultaneously provides a segmentation of $u$. In the new model a Kullback-Leibler divergence term is used to force $v$ to be close to $\bar{v}$, thus allowing it to further extract smaller-scale oscillatory components from the starting cartoon part $\bar{u}$. Thanks to this additional term, the segmentation process is shown to have improved robustness with respect to noise and texture in the initial image.

\noindent The rest of the paper is organized as follows: in Section \ref{sec:CTD} we recall the cartoon-texture decomposition of an image, in Section \ref{sec:CTETRISmodel} we introduce the proposed model, which results in a non-smooth convex optimization problem, and in Section  \ref{sec:minimizingCTETRIS} we introduce an ADMM scheme for the problem solution and analyze its convergence. Section \ref{sec:numericalexp} is devoted to numerical experiments and comparison with the original CEN model and with state-of-the-art models suited for textural image segmentation. Finally, we draw our conclusions in Section \ref{sec:conclusions}.

\section{Cartoon-Texture Decomposition}\label{sec:CTD}
An image $f$ is usually described as a superposition of two components, i.e.,
$$f=u+v,$$
\noindent where $u$ is the geometric component and $v$ is the oscillatory one. The geometric component, commonly referred to as `cartoon', consists of the piecewise-constant parts of an image, including homogeneous regions, contours, and sharp edges.
In contrast, the oscillatory component includes the patterns which can be observed in the image, such as texture or noise. Both texture and noise can indeed be seen as repeated patterns of small scale details, with noise being characterized by random and uncorrelated values. The cartoon-texture decomposition of an image plays an important role in computer vision \cite{bib:CTD}, with a wide range of applications to, e.g., image restoration, segmentation, image editing, and remote sensing. It is an underdetermined linear inverse problem with many solutions, usually described by variational models able to force the cartoon and the texture into different functional spaces in order to produce the required decomposition.

\begin{figure*}[ht!]
\medskip
\begin{center}
\newcolumntype{C}{>{\centering\arraybackslash} m{.25\textwidth} }
\begin{tabular}{CCC}
{\small original image }  &
{\small cartoon }  &
{\small texture }\\
\includegraphics[width=.24\textwidth]{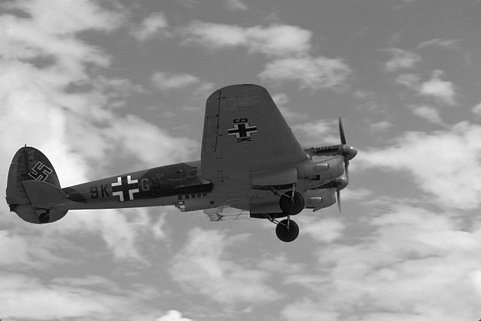} &
\includegraphics[width=.24\textwidth]{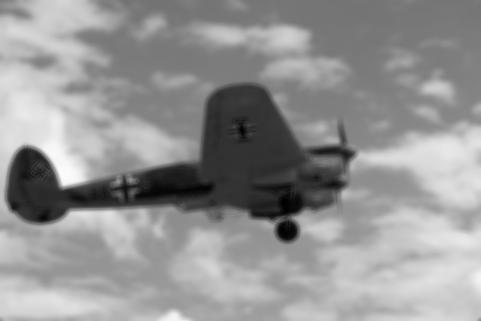} &
\includegraphics[width=.24\textwidth]{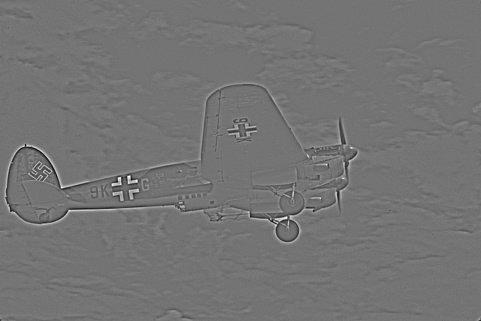}
\end{tabular}
\end{center}
\caption{Cartoon-texture decomposition of \textit{airplane} image after the application of \eqref{ctd_formula}-\eqref{eq:relredrate}.}
\label{fig:ctd}
\end{figure*}

\noindent Following the idea of Meyer \cite{bib:Meyer2001},
 the general image decomposition problem can be formulated as
 \begin{equation}\label{ct_frame}
 \begin{array}{rl}
\underset{(u,v)\in X \times Y} {\min} & g_1(u)+g_2(v)\\
 s.t. & u+v=f,\\
 \end{array}
 \end{equation}
 
\noindent where $X$ and $Y$ are suitable function spaces and $g_1$ and $g_2$ are functionals that model the cartoon regions and the texture patterns, respectively. Several choices have been proposed in literature for both $X, Y$ and  $g_1, g_2$ \cite{bib:CTD1,bib:CTD2}.
A widely used choice to model the cartoon is $g_1(u)=TV(u)$, due to its ability to induce piecewise smooth $u$ with bounded variations \cite{bib:CTD-TV,bib:CTD-TV1}.
Some alternative approaches impose a sparse representation of the cartoon under a given system,
such as wavelet frames~\cite{bib:CTD-SR} or curvelet systems~\cite{bib:CTD-SR1}.
Modeling the texture component is a more complex task, due to the difficulty of conceptualizing mathematical properties able to encompass all the texture types. Many models use the space of oscillatory functions equipped with appropriate norms able to represent
textured or oscillatory patterns \cite{bib:CTD-TV,bib:CTD-TV1,bib:GSpace}. 
An alternative approach assumes that, under suitable conditions, textures can be sparsified, i.e.,
a texture patch can be represented by few atoms in a given dictionary or by specific transforms \cite{bib:Xu2018ImageCD}.\\
Since the existing methods for cartoon-texture decomposition are beyond the scope of this paper, 
here we simply assume that we are able to obtain a decomposition of the given image: 
\begin{equation}\label{ctd}
   \bar{f} = \bar{u} + \bar{v},
\end{equation}
with the aim of using the different information on the two components to improve the effectiveness of the CEN model.
In our experiments we will consider the algorithm described in \cite{bib:cartoontexture}. Figure \ref{fig:ctd} shows the decomposition produced by one iteration of the algorithm, which results in 
\begin{equation}\label{ctd_formula}
  \bar{u}(x)=\omega(\rho_\sigma(x))L_\sigma * \bar{f}+(1-\omega(\rho_\sigma(x)))\bar{f}, \;\;\;\; \bar{v}(x)=\bar{f}(x)-\bar{u}(x),  
\end{equation}

\noindent where  $L_\sigma$  is a low-pass filter, $*$ is the convolution operator, $\omega: [0,1] \longrightarrow [0,1]$
 is an increasing function that is constant and equal to zero
near zero and constant and equal to 1 near 1, and $\rho_\sigma(x)$ is the relative reduction rate of local TV
\begin{equation}\label{eq:relredrate}
\rho_\sigma(x) =\frac{LTV_\sigma (\bar{f}(x))-LTV_\sigma (L_\sigma *\bar{f}(x))}{LTV_\sigma (\bar{f}(x))} \in [0,1]
\end{equation}

\noindent with
$LTV_\sigma (\bar{f}(x)) = \left( L_\sigma * \vert \nabla \bar{f}\vert \right)$.

\noindent We note that the cartoon-texture decomposition produced by \eqref{ctd_formula} is not unique, but it depends on the choice of $\sigma$ \cite{bib:cartoontexture}. Anyway, we will show that a rough decomposition is enough for our model, hence there's no need for an accurate tuning of $\sigma$.

\section{The C-TETRIS model}\label{sec:CTETRISmodel}
We here introduce the Cartoon-Texture Evolution for Two-Region Image Segmentation (C-TETRIS) model. As mentioned in the previous sections, starting from the decomposition \eqref{ctd}, the main idea behind C-TETRIS is to simultaneously produce the segmentation of $\bar{u}$ and its cartoon-texture decomposition. 
In detail, it decomposes $\bar{u}$ as $\bar{u}=u+v$, where $v$ is enforced to be close to $\bar{v}$, and computes a segmentation of $u$ by solving the problem 
\begin{equation}\label{ctetris}
\begin{array}{rl}
\underset{u,c_{in},c_{out},v}{\min} & \displaystyle {\mathcal E}_{CEN}(u,c_{in},c_{out}; \bar{u}) + \mu {\mathcal D}_{KL}(v;\bar{v})\\
\stt     & 0 \leq u \leq  1,\\
& \displaystyle u+v=\bar{u},
\end{array}
\end{equation}

\noindent where ${\mathcal E}_{CEN}$ represents the objective function of problem \eqref{eq:cen_continuous}, ${\mathcal D}_{KL}(v;\bar{v})$ denotes the Kullback-Leibler (KL) divergence of $v$ from $\bar{v}$, defined as 
\begin{equation}\label{KL}
    {\mathcal D}_{KL}(v;\bar{v}) =\int_{\Omega} v(x) \log \left(\frac{v(x)}{\bar{v}(x)}\right) dx,
\end{equation}	
\noindent where we set
\[v(x) \log \left(\frac{v(x)}{\bar{v}(x)} \right) = \left \{ \begin{array}{cc}
 0    & v(x)=0, \\
\infty     &  \bar{v}(x) = 0,
\end{array} 
\right.
\]
\noindent and $\mu>0$. The KL divergence measures the amount of information lost if $\bar{v}$ is used to approximate $v$ and appears in many models of imaging science, where it is usually employed as a fidelity term.
Simply speaking, the C-TETRIS model extracts from $\bar{u}$ the ``remaining texture'' and produces its best approximation among all the functions that take only two values. \\
In the following we consider the discrete version of (\ref{ctetris}).
 Let
$$
    \Omega_{n_x,n_y}=\left\{ (i,j) : 0 \leq i \leq n_x-1, \, 0 \leq j \leq n_y-1\right\}
$$
be a discretization of $\Omega$ consisting of an $ n_x \times n_y$ grid of pixels
and
$$
    \vert \nabla_x u \vert_{i,j}  = \vert \delta_x^+ u \vert_{i,j} , \quad  \vert \nabla_y u \vert_{i,j} = \vert \delta_y^+ u \vert_{i,j}
$$ 
where $\delta_x^+$  and $\delta_y^+$ are the forward finite-difference operators in the $x$- and $y$-directions, with unit spacing,
and the values $u_{i,j}$ with indices outside~$\Omega_{n_x,n_y}$ are defined by replication. 
The discrete version of the (\ref{ctetris}) leads to the following non-smooth constrained optimization problem:

\begin{equation} \label{eq:ctetris_discr}
\begin{array}{rl}
    \underset{u, c_{in}, c_{out},v}{\min} &  \displaystyle E_{CEN}(u,c_{in},c_{out}; \bar{u}) + \mu D_{KL}(v;\bar{v}) \\
   \stt             & 0 \le u \le 1, \\[2pt]
                    & u +v =\bar u,
\end{array}
\end{equation}

\noindent where we denoted by $E_{CEN}$ the discrete version of ${\mathcal E}_{CEN}$, defined as 
\begin{multline*}
    E_{CEN}(u,c_{in},c_{out}; \bar{u}) = \sum_{i,j} \big(\vert \nabla_x u \vert_{i,j} + \vert \nabla_y u \vert_{i,j}\big) +\\ +\lambda \sum_{i,j} \left( u_{i,j} (c_{in}-\bar{u}_{i,j})^2 + (1- u_{i,j})\,( c_{out}-\bar{u}_{i,j})^2\right),
\end{multline*}

\noindent and we denoted with $D_{KL}$ the discrete version of the Kullback-Leibler divergence ${\mathcal D}_{KL}$, defined as
$$D_{KL}(v;\bar{v}) = \sum_{i,j} v_{i,j} \log \left( \frac{v_{i,j}}{\bar{v}_{i,j}}\right).$$

\noindent It is worth noting that the first term in $E_{CEN}$ corresponds to the discrete Total Variation (TV) of the image $u$. We here opted for the use of a modified version of the TV functional, in which the $\ell_2$ norm is replaced by the $\ell_1$ one (as proposed in \cite{bib:Esedoglu2004}), 
since in the case of image restoration it is known to produce sharper piece-wise constant images. 
Nevertheless, a preliminary comparison between the models equipped with the $\ell_1$ and the $\ell_2$ version, respectively, showed no difference in terms of segmentation quality.

\section{Minimizing the C-TETRIS model}\label{sec:minimizingCTETRIS}

We here focus on the solution of the minimization problem in \eqref{eq:ctetris_discr}. One can observe that, although the problem is in general nonconvex, it becomes convex when either the pair $(c_{in}, \,c_{out})$ or the pair $(u,\,v)$ are fixed. Suppose, for the moment, that the values of $c_{in}, c_{out}$ have been determined and consider the minimization problem in $u$ and $v$ only, which can be written as
\begin{equation} \label{eq:ctetris_discr_Cfixed}
\begin{array}{rl}
    \underset{u,v}{\min} &  \displaystyle  \sum_{i,j} \big(\vert \nabla_x u \vert_{i,j} + \vert \nabla_y u \vert_{i,j}\big) +
                         \lambda \, r^\top  u + \mu D_{KL}(v;\bar{v}) \\
   \stt             & 0 \le u \le 1, \\[2pt]
                    & u + v =\bar u,
\end{array}
\end{equation}

\noindent where we defined, for each $(i,j)$,
$$r_{i,j}\equiv r_{i,j}(c_{in},c_{out})=    \left( c_{in}-\bar{u}_{i,j} \right)^2 - \left( c_{out}-\bar{u}_{i,j}\right)^2.$$

\noindent Problem \eqref{eq:ctetris_discr_Cfixed} is a non-smooth convex optimization problem subject to linear and bound constraints which we propose to solve by the Alternate Directions Method of Multipliers (ADMM) \cite{bib:BoydEtAl2011admm}. To this aim, we reformulate problem \eqref{eq:ctetris_discr_Cfixed} as
\begin{equation} \label{eq:ctetris_discr_Cfixed_ref1}
\begin{array}{rl}
    \underset{u,d_x,d_y,v}{\min} &  \displaystyle \|d_x\|_1 + \|d_y\|_1 + \lambda \, r^\top  u  + \mu D_{KL}(v;\bar{v}) \\
   \stt          & d_x = \nabla_x u,   \\[2pt]
                & d_y = \nabla_y u,   \\[2pt]
                & u + v =\bar u, \\[2pt]
                & 0 \le u \le 1.
\end{array}
\end{equation}

\noindent Starting from \eqref{eq:ctetris_discr_Cfixed_ref1}, it is straightforward to check that the objective function and the constraints of the problem can be split in two blocks. Indeed, by introducing the variable $z = [d_x^\top,d_y^\top,v^\top]^\top$, one can further reformulate  \eqref{eq:ctetris_discr_Cfixed_ref1} as

\begin{equation} \label{eq:ctetris_discr_Cfixed_ref2}
\begin{array}{rl}
    \underset{u,z}{\min} &  \displaystyle F(u) + G(z) \\
   \stt & H\,u - z = b,
\end{array}
\end{equation}
where we defined
\begin{eqnarray*}
&F(u) = \lambda \, r^\top  u + \chi_{[0,1]}(u),\qquad G(z) = \|d_x\|_1 + \|d_y\|_1 + \mu D_{KL}(v;\bar{v}),\\
&H = \left[ \nabla_x^\top,\,\nabla_y^\top,\,-I \right]^\top, \mbox{and }\quad b = [ 0,\,0,\,-\bar{u}^\top]^\top,
\end{eqnarray*}
and we used $\chi_{[0,1]}(u)$ to indicate the characteristic function of the hypercube $[0,1]^{n_x\times n_y}$.

\noindent Consider the Lagrangian and the augmented Lagrangian functions associated with problem \eqref{eq:ctetris_discr_Cfixed_ref2}, defined respectively as
\begin{eqnarray*}
    &\mathcal{L}(u,z,\xi) = F(u) + G(z) + \xi^\top \left(H\,u - z - b\right),\\[2mm]
    &\mathcal{L}_A(u,z,\xi;\rho) = F(u) + G(z) + \xi^\top \left(H\,u - z - b\right) + \frac{\rho}{2}\left\|H\,u - z - b\right\|_2^2,
\end{eqnarray*}
\noindent
where $\rho>0$, and $\xi$ is a vector of Lagrange multipliers.

\noindent Starting from given estimates $u^0$, $z^0$, and $\xi^0$, at each iteration $k$ ADMM updates the estimates as
\begin{equation}\label{eq:admm_method}
\begin{split}
    u^{k+1} & = \displaystyle \argmin\limits_{u} \mathcal{L}_A(u,z^k,\xi^k;\rho),\\
    z^{k+1} & = \displaystyle \argmin\limits_{z} \mathcal{L}_A(u^{k+1},z,\xi^k;\rho),\\
    \xi^{k+1} & = \displaystyle \xi^k + \rho\left(H\,u^{k+1} - z^{k+1}\right).
\end{split}
\end{equation}

\noindent Since $F(u)$ and $G(z)$ in~\eqref{eq:ctetris_discr_Cfixed_ref2} are closed, proper and convex, and $H$ has full rank, the convergence of ADMM can be proved by exploiting the classical result from \cite{bib:Eckstein1992}, which we report in the following.

\begin{theorem}\label{thm:ADMM_conv}
Consider problem \eqref{eq:ctetris_discr_Cfixed_ref2} where $F(u)$ and $G(z)$ are closed, proper
and convex functions and $H$ has full rank. Consider the summable sequences $\{\varepsilon_k\}, \{\nu_k\} \subset \R_+$ and let

\begin{eqnarray*}
&&\left\| u^{k+1} - \argmin\limits_{u} \mathcal{L}_A(u,z^k,\xi^k;\rho)\right\| \leq \varepsilon_k,\\
&&\left\| z^{k+1} - \argmin\limits_{z} \mathcal{L}_A(u^{k+1},z,\xi^k;\rho)\right\| \leq \nu_k,\\
&&\xi^{k+1} = \xi^k + \rho\left(H\,u^{k+1} - z^{k+1}\right).
\end{eqnarray*}

\noindent
If there exists a saddle point $(u^*,z^*,\xi^*)$ of $ \mathcal{L}(u,z,\xi)$, then $u^k\rightarrow u^*$,
$z^k\rightarrow z^*$ and $\xi^k\rightarrow\xi^*$. If such saddle point does not exist, then at least
one of the sequences $\{z^k\}$ or $\{\xi^k\}$ is unbounded.
\end{theorem}

\noindent Theorem~\ref{thm:ADMM_conv} guarantees the convergence of the ADMM scheme even if the subproblems are solved inexactly, provided that the inexactness of the solution can be controlled.

\noindent So far we have been concerned with the solution of problem \eqref{eq:ctetris_discr_Cfixed} when the values of $c_{in}$ and $c_{out}$ are known in advance which, however, is not the case in practice. By following the example of \cite{bib:CEN2006}, we adopt a two-step scheme in which we alternate updates of $u$ and $z$, determining the shape of the two regions, and updates of $c_{in}$ and $c_{out}$. Observe that, by fixing $u=u^k$ and $z=z^k$, the restriction of problem \eqref{eq:ctetris_discr} to $c_{in}$ and $c_{out}$ can be written as the unconstrained convex quadratic optimization problem

\begin{equation} \label{eq:ctetris_discr_fixUV}
    \underset{c_{in}, c_{out}}{\min} \quad  \displaystyle \sum_{i,j} \left( u^k_{i,j} (c_{in}-\bar{u}_{i,j})^2 + (1- u^k_{i,j})\,( c_{out}-\bar{u}_{i,j})^2\right).
\end{equation}
Hence, we propose to update the values of $c_{in}$ and $c_{out}$ after each ADMM step by taking the exact minimizer of problem \eqref{eq:ctetris_discr_fixUV}, i.e., by setting

\begin{equation}\label{eq:updateC}
    c_{in}^k = \frac{\sum_{i,j}u^k_{i,j}\bar{u}_{i,j}}{\sum_{i,j}u^k_{i,j}},\quad\mbox{and}\quad c_{out}^k = \frac{\sum_{i,j}(1-u^k_{i,j})\bar{u}_{i,j}}{\sum_{i,j}(1-u^k_{i,j})}.    
\end{equation}

It is worth pointing out that such a modification alters the original ADMM scheme making it an inexact alternate minimization scheme for the problem in $u$, $z$, $c_{in}$, and $c_{out}$. Nevertheless, as also shown for the original CEN model, the experiments carried out in this work show that in all the cases under analysis the values of $c_{in}$ and $c_{out}$ stagnate after the first few iterations, thus recovering in practice the convergence properties shown for the case of fixed $c_{in}$ and $c_{out}$.

\subsection{Solving the ADMM subproblems}
We will now focus on how the subproblems in \eqref{eq:admm_method} can be solved in practice. First, by expliciting the form of the augmented Lagrangian functions, we can rewrite the ADMM scheme as

\begin{equation*}
\begin{split}
    u^{k+1} & = \displaystyle \argmin\limits_{0\leq u \leq 1} \lambda \, r^\top  u + (\xi^k)^\top \left(H\,u - z^k - b\right) + \frac{\rho}{2}\left\|H\,u - z^k - b\right\|_2^2,\\
    z^{k+1} & = \displaystyle \argmin\limits_{z} G(z) + (\xi^k)^\top \left(H\,u^{k+1} - z - b\right) + \frac{\rho}{2}\left\|H\,u^{k+1} - z - b\right\|_2^2,\\
    \xi^{k+1} & = \displaystyle \xi^k + \rho\left(H\,u^{k+1} - z^{k+1} - b\right).
\end{split}
\end{equation*}

\noindent It is straightforward to check that the minimization problem over $z$ can be split into three independent minimization problems, respectively on $d_x$, $d_y$, and $v$, leading to the following scheme

\begin{equation}\label{eq:admm_method_4blocks}
\begin{split}
    u^{k+1} & = \displaystyle \argmin\limits_{0\leq u \leq 1} \lambda \, r^\top  u + (\xi^k)^\top \left(H\,u - z^k - b\right) + \frac{\rho}{2}\left\|H\,u - z^k\right\|_2^2,\\
    d_x^{k+1} & = \displaystyle \argmin\limits_{d_x} \|d_x\|_1 + (\xi_x^k)^\top \left(\nabla_x u^{k+1} - d_x\right) + \frac{\rho}{2}\left\|\nabla_x u^{k+1} - d_x\right\|_2^2,\\
    d_y^{k+1} & = \displaystyle \argmin\limits_{d_y} \|d_y\|_1 + (\xi_y^k)^\top \left(\nabla_y u^{k+1} - d_y\right) + \frac{\rho}{2}\left\|\nabla_y u^{k+1} - d_y\right\|_2^2,\\
    v^{k+1} & = \displaystyle \argmin\limits_{v} \mu D_{KL}(v;\bar{v}) + (\xi_v^k)^\top \left(-u^{k+1} - v + \bar{u}\right) +  \frac{\rho}{2}\left\|u^{k+1} + v - \bar{u}\right\|_2^2,\\
    \xi^{k+1} & = \displaystyle \xi^k + \rho\left(H\,u^{k+1} - z^{k+1} -b\right),
\end{split}
\end{equation}

\noindent where we split the Lagrange multipliers vector $\xi$ as $\xi = [\xi_x^\top, \xi_y^\top, \xi_v^\top]^\top $. The scheme presented in \eqref{eq:admm_method_4blocks} can be further simplified by exploiting the linearity of the constraints $H\,u -z = b$, as suggested in \cite{bib:GoldsteinOsher2009}. In detail, by introducing the vectors $ b_x^k = \frac{\xi_x^k}{\rho}$, $b_y^k = \frac{\xi_y^k}{\rho}$, and $ b_v^k = -\frac{\xi_v^k}{\rho}-\bar{u}$, one can rewrite \eqref{eq:admm_method_4blocks} equivalently as

\begin{align}
\begin{split}\label{eq:admm_method_4blocks_simpU}
        u^{k+1} =&\;\; \displaystyle \argmin\limits_{0\leq u \leq 1} \lambda \, r^\top  u + \frac{\rho}{2}\left\|\nabla_x u - d_x^k + b_x^k\right\|_2^2 \\
        &\qquad\qquad+ \frac{\rho}{2}\left\|\nabla_y u - d_y^k + b_y^k\right\|_2^2 + \frac{\rho}{2}\left\| u + v^k + b_v^k\right\|_2^2,
\end{split}\\
    d_x^{k+1} =&\;\; \displaystyle \argmin\limits_{d_x}  \|d_x\|_1 + \frac{\rho}{2}\left\|\nabla_x u^{k+1} - d_x + b_x^k\right\|_2^2, \label{eq:admm_method_4blocks_simpDX}\\
    d_y^{k+1} =&\;\; \displaystyle \argmin\limits_{d_y}  \|d_y\|_1 + \frac{\rho}{2}\left\|\nabla_y u^{k+1} - d_y + b_y^k\right\|_2^2, \label{eq:admm_method_4blocks_simpDY}\\
    v^{k+1} =&\;\; \displaystyle \argmin\limits_{v}  \mu D_{KL}(v;\bar{v}) + \frac{\rho}{2}\left\|u^{k+1} + v + b_v^k\right\|_2^2, \label{eq:admm_method_4blocks_simpV}\\[2mm]
    \begin{split}\label{eq:admm_method_4blocks_simpB}
        b_x^{k+1} =&\;\;\; \displaystyle b_x^k + \nabla_x\,u^{k+1} - d_x^{k+1},\\[2mm]
        b_y^{k+1} =&\;\;\; \displaystyle b_y^k + \nabla_y\,u^{k+1} - d_y^{k+1},\\[2mm]
        b_v^{k+1} =&\;\;\; \displaystyle b_v^k + u^{k+1} + v^{k+1} - \bar{u}.
    \end{split}
\end{align}

\noindent Problem \eqref{eq:admm_method_4blocks_simpU} is a strongly convex bound-constrained quadratic optimization problem. To obtain an approximate solution $u^{k+1}$, by following \cite{bib:GoldsteinBressonOsher2010,bib:Antonelli2020Adaptive}, we consider the optimality conditions of the unconstrained version of the problem, i.e., the solution to the linear system
$$(- \Delta + I) u  = - \frac{\lambda\,r}{\rho} + ( \nabla_x^\top ( b_x^{k} - d_x^{k}) )+ ( \nabla_y^\top ( b_y^{k} - d_y^{k}) ) + (b^{k}-v^{k}),$$
where $\Delta$ represents the finite-difference discretization of the Laplacian. We first solve the system by Gauss-Seidel method and then project the solution in $[0,1]^{n_x \times n_y}$.

\noindent As regards the updates in \eqref{eq:admm_method_4blocks_simpDX}-\eqref{eq:admm_method_4blocks_simpV}, one has to note that they are proximal operators \cite{bib2014:ParikhProxAlg, bib:Beck2017ProxAlg} of closed proper and convex functions. In detail, the proximal operator in \eqref{eq:admm_method_4blocks_simpDX} and \eqref{eq:admm_method_4blocks_simpDY} can be computed in closed form by means of the well-known soft-thresholding operator, defined as
$$
[{\mathcal S}(x,\gamma)]_{i,j}= \mathrm{sign}(x_{i,j})\cdot\max\big(\vert x_{i,j}\vert-\gamma, 0\big).
$$
Finally, the proximal operator in \eqref{eq:admm_method_4blocks_simpV} can be computed as
$$ [\mathrm{prox}_{\gamma D_{KL}(x,\tilde{x})} (x)]_{i,j} = \gamma 
 W (\gamma^{-1}\tilde{x}_{i,j} e^{\gamma^{-1}x_{i,j}-\tilde{x}_{i,j}^{-1}}), $$
\noindent where $W (x)$ is the Lambert $W$ function satisfying $W (y) e^{W(y)} = y$ which, although not available in closed form, can be approximated with high precision.

\section{Numerical experiments}\label{sec:numericalexp}
\begin{figure}[!t]
\medskip
\begin{center}
\newcolumntype{C}{>{\centering\arraybackslash} m{.20\textwidth} }
\begin{tabular}{C|CCC}
 cartoon-texture decomposition & \multicolumn{3}{c}{ ADMM iterations} \\
 & {\small first} & {\small intermediate} & {\small last}\\
 \hline\\
 {\small cartoon} & & & \\
 \includegraphics[width=.17\textwidth]{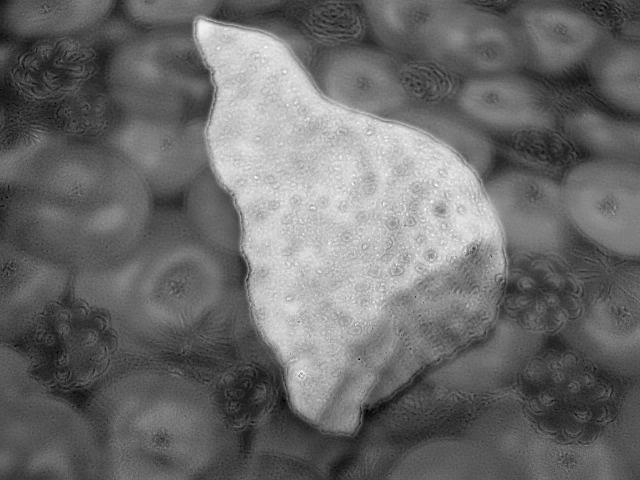} &
\includegraphics[width=.17\textwidth]{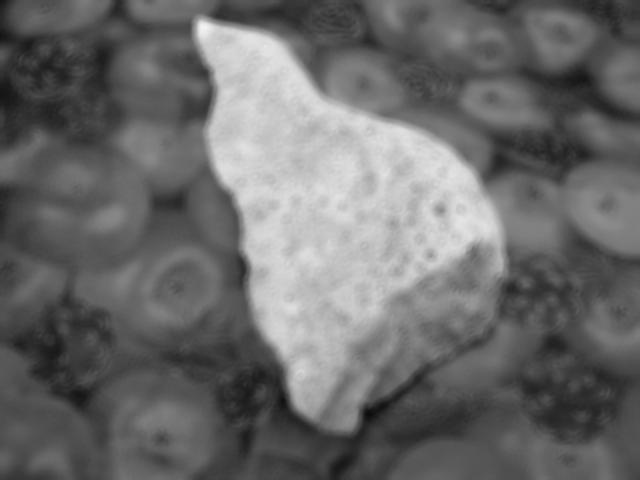} &
\includegraphics[width=.17\textwidth]{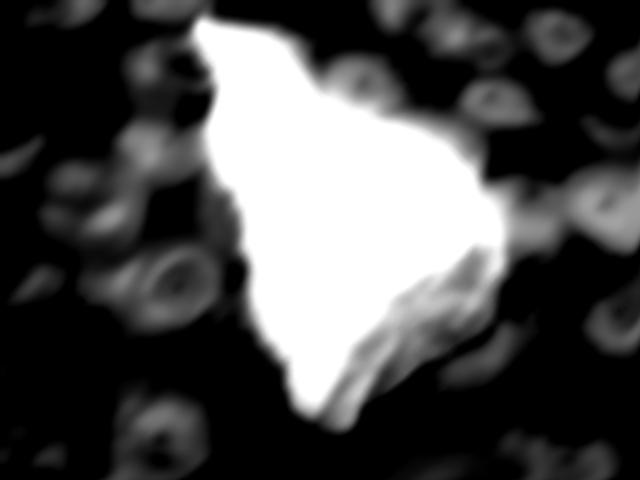} &
\includegraphics[width=.17\textwidth]{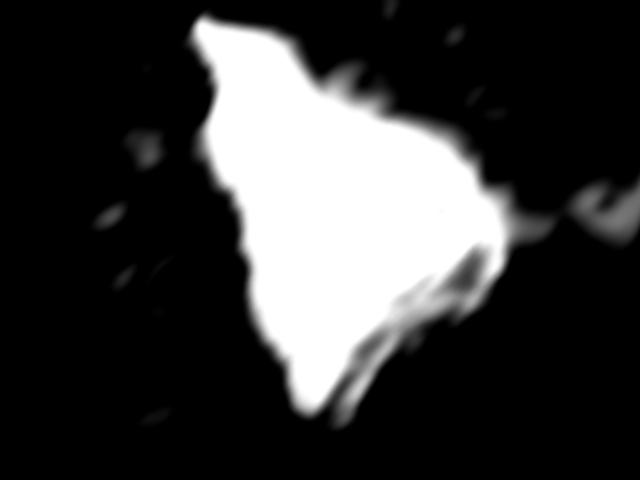}  \\
{\small texture} & & & \\
 \includegraphics[width=.17\textwidth]{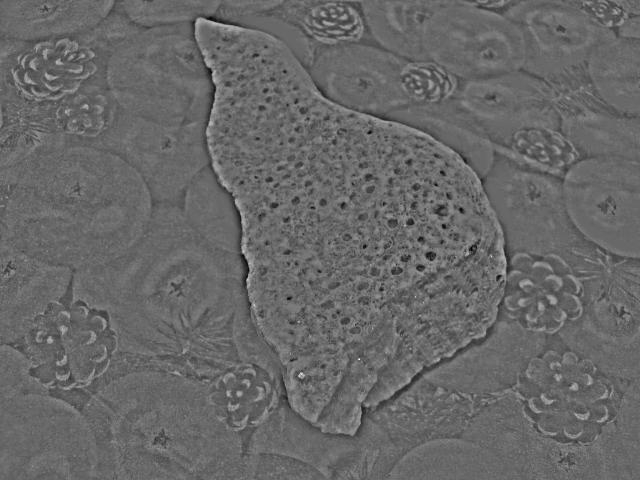} &
\includegraphics[width=.17\textwidth]{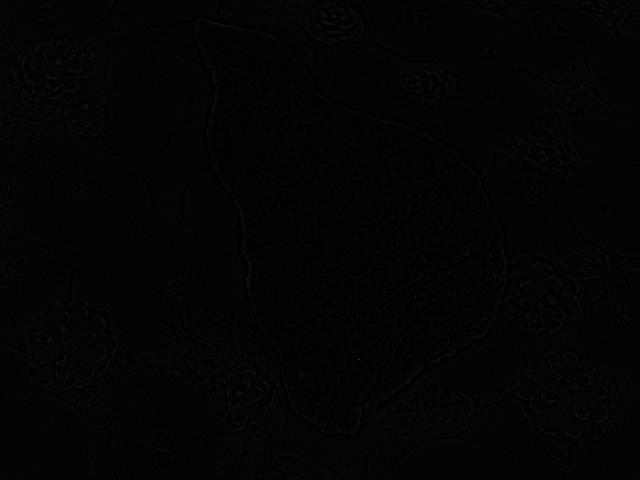} &
\includegraphics[width=.17\textwidth]{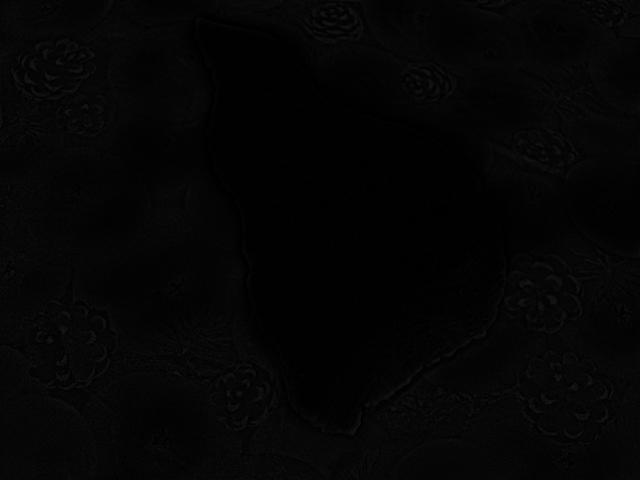} &
\includegraphics[width=.17\textwidth]{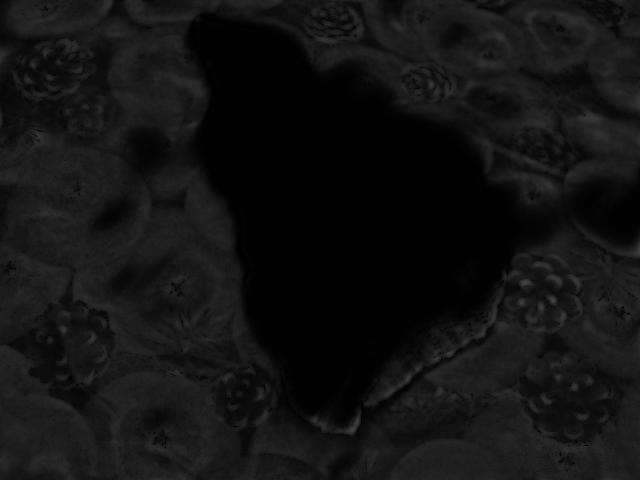}  \\
\hline 
 {\small cartoon} & & & \\
 \includegraphics[width=.17\textwidth]{airplane2U_SBit0.jpg} &
\includegraphics[width=.17\textwidth]{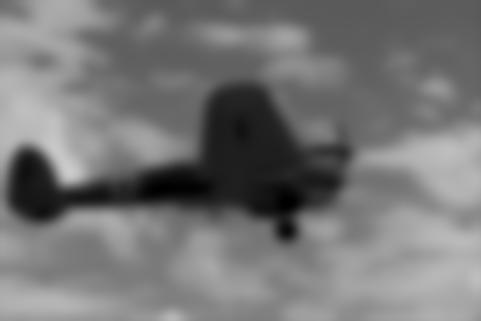} &
\includegraphics[width=.17\textwidth]{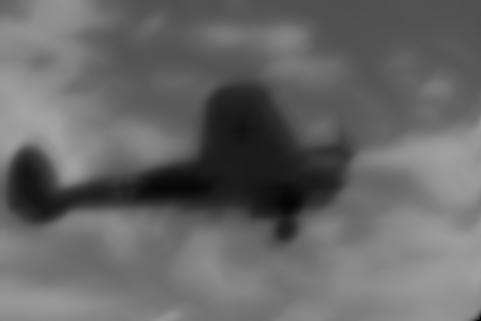} &
\includegraphics[width=.17\textwidth]{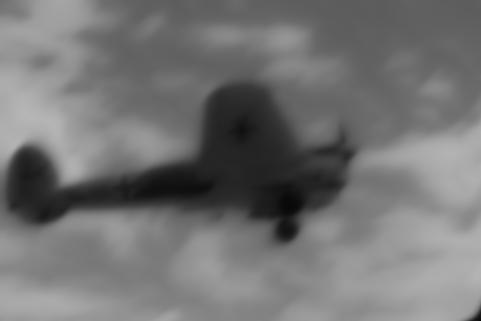}  \\
{\small texture} & & & \\
 \includegraphics[width=.17\textwidth]{airplane2V_SBit0.jpg} &
\includegraphics[width=.17\textwidth]{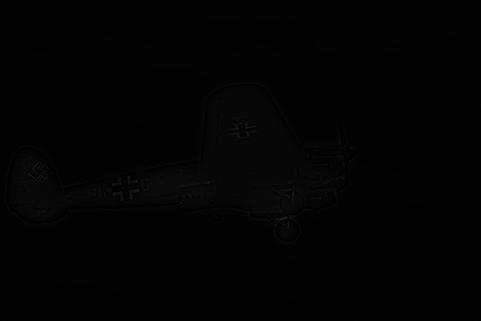} &
\includegraphics[width=.17\textwidth]{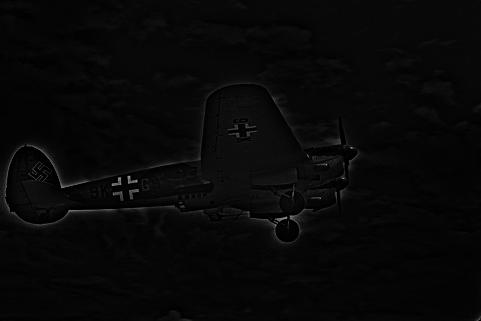} &
\includegraphics[width=.17\textwidth]{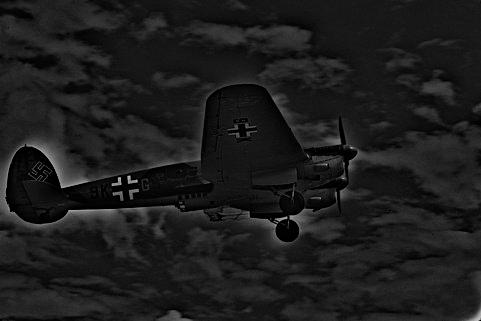} \\
\hline 
{\small cartoon} & & & \\
 \includegraphics[width=.17\textwidth]{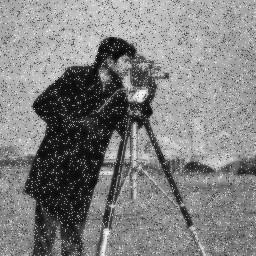} &
\includegraphics[width=.17\textwidth]{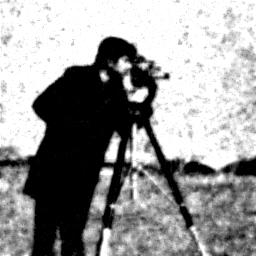} &
\includegraphics[width=.17\textwidth]{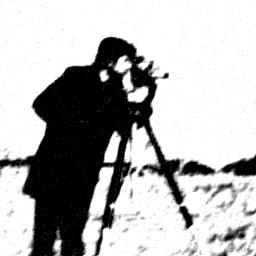} &
\includegraphics[width=.17\textwidth]{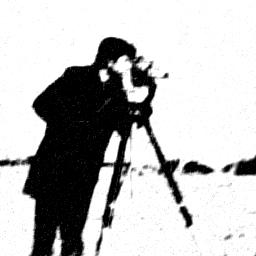}  \\
{\small texture} & & & \\
 \includegraphics[width=.17\textwidth]{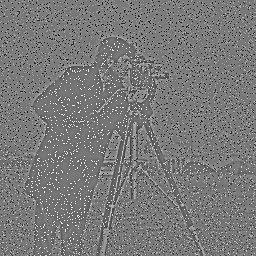} &
\includegraphics[width=.17\textwidth]{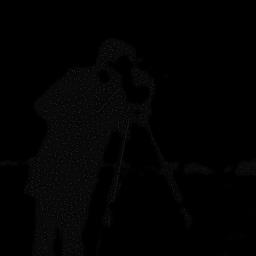} &
\includegraphics[width=.17\textwidth]{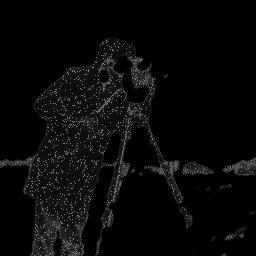} &
\includegraphics[width=.17\textwidth]{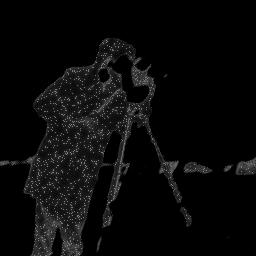} \\
\hline
{\small cartoon} & & & \\
\includegraphics[width=.17\textwidth]{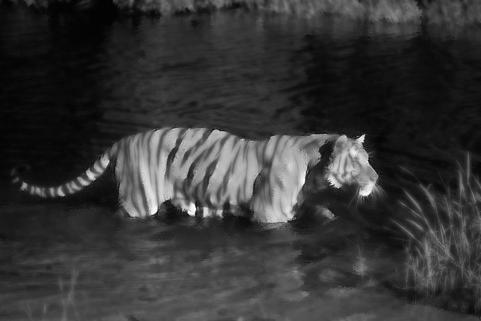} &
\includegraphics[width=.17\textwidth]{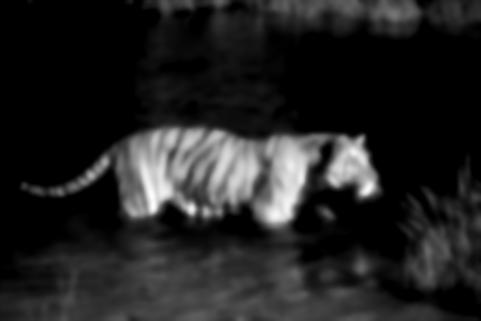} &
\includegraphics[width=.17\textwidth]{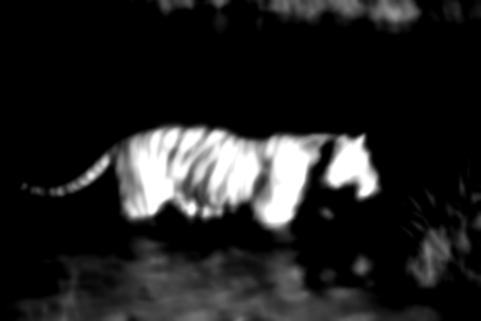} &
\includegraphics[width=.17\textwidth]{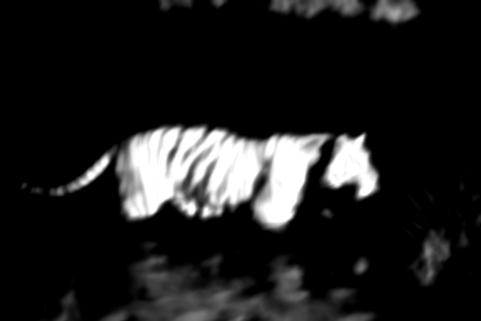} \\
{\small texture} & & & \\
\includegraphics[width=.17\textwidth]{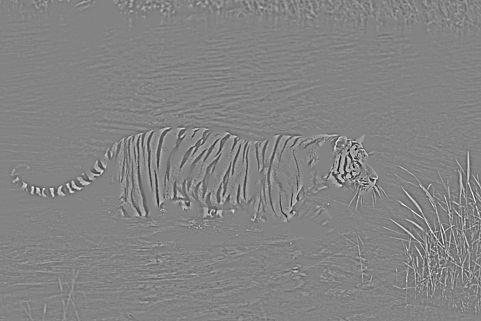} &
\includegraphics[width=.17\textwidth]{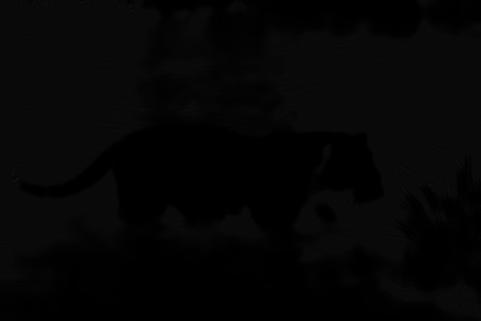} &
\includegraphics[width=.17\textwidth]{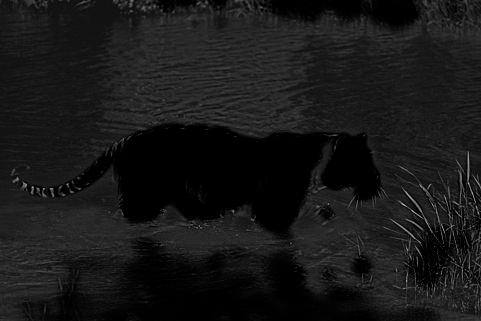} &
\includegraphics[width=.17\textwidth]{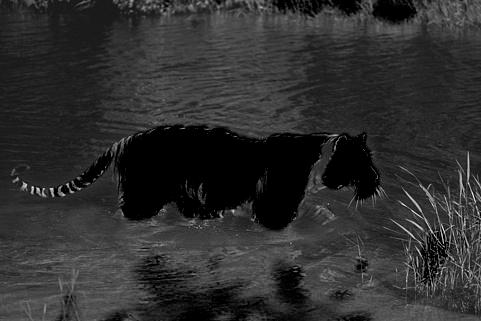} \\
\end{tabular}
\vskip 1pt
\caption{Details of the evolution of the cartoon ($u$) and the texture part ($v$) performed by C-TETRIS on images \textit{airplane}, \textit{cameraman} with 15$\%$ of salt \& pepper noise, and \textit{tiger}. The segmentations produced on the last iteration for each image are showed in Figg. \ref{fig:gt}, \ref{fig:CENvsCTETRIS}, \ref{fig:noise}, and \ref{fig:texture}, respectively.
\label{fig:CTevolution}}
\end{center}
\end{figure}

In this section, we test the effectiveness of C-TETRIS in producing two-region segmentation on various image sets. The first set contains three pairs of real-life images with corresponding ground truth coming from the database \cite{bib:datasetGT}: \textit{man} is a smooth image whereas \textit{flowerbed} and \textit{stone} show an object foreground on a textured background. The second set consists of four images (see Figure~\ref{fig:CENvsCTETRIS}) available from the Berkeley database \cite{bib:Arbelaez2011BK500} which are in general considered to be smooth: the real-life images \textit{airplane} and \textit{squirrel}, and the medical images \textit{brain} and \textit{ultrasound}. The third set of images consists of noisy versions of the famous \textit{cameraman} image from MIT Image Library\footnote{https://libguides.mit.edu/findingimages} (see Figure~\ref{fig:noise}) which we use to test the robustness of the C-TETRIS model with respect to the noise. The fourth and last set of images (see Figure~\ref{fig:texture}) consists of three textural images: \textit{tiger} and \textit{bear}, taken from  \cite{bib:Arbelaez2011BK500}, and \textit{spiral}, taken from \cite{bib:ChanVese2001}. 
We here provide some further details on the numerical experiments.
The C-TETRIS algorithm was implemented in MATLAB using the Image Processing Toolbox, where the cartoon-texture decomposition was initially performed by one iteration of the algorithm described in~\cite{bib:cartoontexture}, using a Gaussian filter with $\sigma = 2$ as $L_{\sigma}$, and the following function $\omega$ \cite{bib:cartoontexture}:
\begin{equation}
    \omega(x) = \left\{ \begin{array}{ll}
         0,                  & x \leq l_1, \\
         (x-l_1)/(l_2 -l_1), & l_1 \lt x \lt l_2, \\
         1,                  & x \geq l_2,
    \end{array}\right. 
\end{equation}
where the weights $l_1$ and $l_2$ have been set to $0.25$ and $0.5$, respectively. We would like to remark that extensive testing showed that the accuracy of the produced segmentation is only slightly influenced by the variation of the Gaussian smoothing parameter, $\sigma$, or by the number of steps performed to obtain the cartoon-texture decomposition. 
\noindent Among the several available implementations of CEN we chose the one\footnote{\label{footnotebresson}\url{http://htmlpreview.github.io/?https://github.com/xbresson/old_codes/blob/master/codes.html}} proposed by the authors of \cite{bib:GoldsteinBressonOsher2010}. Although the code is written in C programming language, a MEX interface is available for testing in MATLAB.
This implementation is based on split Bregman iterations with the following stopping criterion:
\begin{equation}\label{eq:stop-crit}
\vert \mathtt{diff}^k - \mathtt{diff}^{k-1} \vert \le \mathtt{tol} \quad  \mbox{and} \quad k > \mathtt{maxit} ,
\end{equation}
where
$$
     \mathtt{diff}^k = \frac{\mathtt{sd}(f^k)}{\mathtt{sd}(f^k) \cdot \mathtt{sd}(f^{k-1})}, \quad
     \mathtt{sd}(f^k) = \sum_{i,j} (f_{i,j}^k - f_{i,j}^{(k-1)})^2,
$$
$\mathtt{tol}$ is a given tolerance and $\mathtt{maxit}$ is the maximum number of SB iterations. In order to make a fair comparison, all the algorithms presented in the next section use the stopping criterion (\ref{eq:stop-crit}), where we set $\mathtt{maxit}=50$ and $\mathtt{tol} =10^{-6}$ ($\mathtt{tol} =10^{-8}$ for the noisy images). The parameter $\lambda$ in \eqref{eq:cen_continuous} and in \eqref{eq:ctetris_discr_Cfixed}, has a scaling role and was set according to the level of required details in the segmentation. In particular, in each test for CEN model we used the value proposed by the authors in the available code, which we indicate as $\lambda_{CEN}$, based on this empirical rule: $\lambda_{CEN} =10^a$ with $a \in \{-1,0,1\}$ from larger to smaller regularization/smoothing. To balance the presence of the KL term, for C-TETRIS we perform a grid search and select a parameter $\lambda$ with a variation of at most $5\%$ from $\lambda_{CEN}$.
The parameter $\mu$ was set as $\mu =10^c$ with $c \in \{-2,-1,0\}$. Finally, the Bregman parameter $\rho$ was set to $1$. \\
Before proceeding with the experiments on the four image sets described above, we show an example of the functioning of the proposed model. We consider an image for each of the four sets and report in Figure~\ref{fig:CTevolution} the starting cartoon-texture decomposition and the components $u$ and $v$ after the first ADMM iteration, at an intermediate iteration and at the last iteration. We note that, as the ADMM advances, the remaining texture is progressively subtracted from the cartoon, allowing a clearer distinction of background and foreground.

\subsection{Results on ground truth images\label{subsec:groundtruth}}
\begin{figure*}[!ht]
\medskip
\begin{center}
\newcolumntype{C}{>{\centering\arraybackslash} m{.20\textwidth} }
\begin{tabular}{CCCC}
\multicolumn{1}{c}{\small original image }  &
\multicolumn{1}{c}{\small ground truth }   &
\multicolumn{1}{c}{\small CEN } &
\multicolumn{1}{c}{\small C-TETRIS }\\
\hline\\
\includegraphics[width=.18\textwidth]{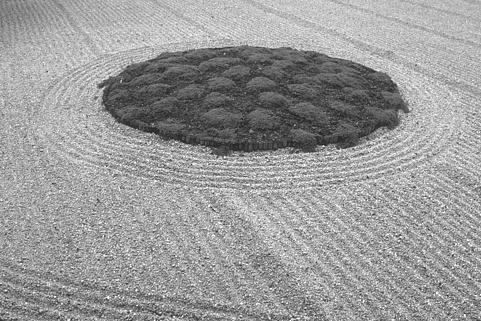} &
\includegraphics[width=.18\textwidth]{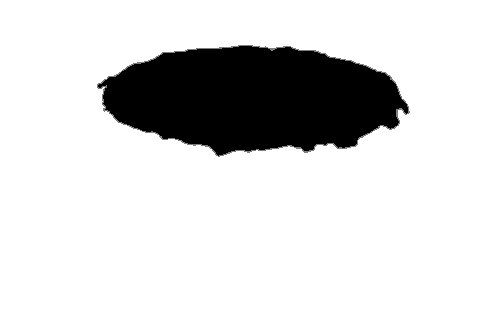}  & \includegraphics[width=.18\textwidth]{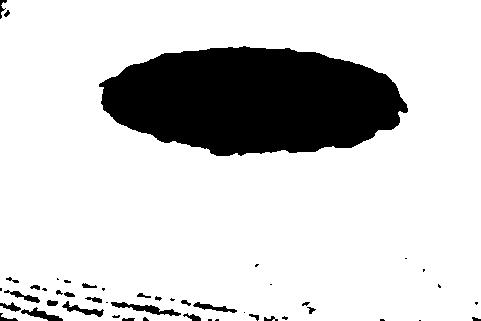} &
\includegraphics[width=.18\textwidth]{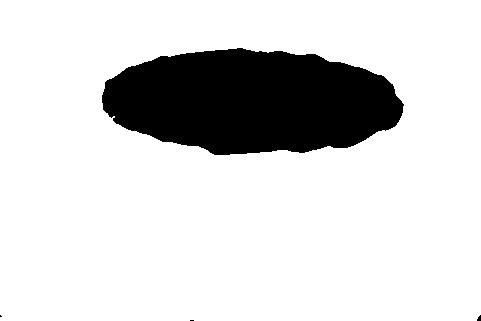} \\[1mm]
\multicolumn{1}{c}{\scriptsize \textit{flowerbed}} & & &\\

\hline\\
\includegraphics[width=.18\textwidth]{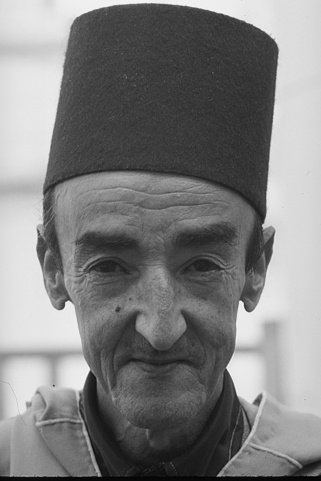} &
\includegraphics[width=.18\textwidth]{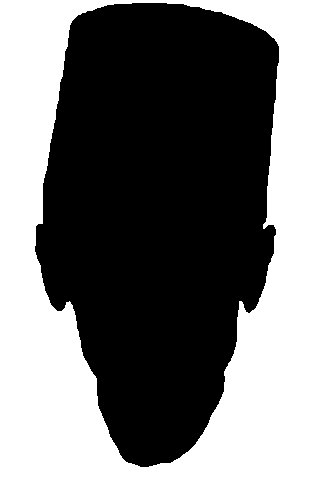}  & \includegraphics[width=.18\textwidth]{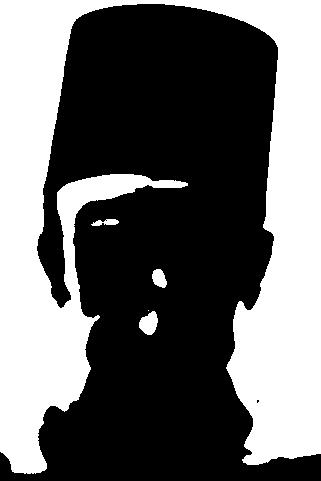} &
\includegraphics[width=.18\textwidth]{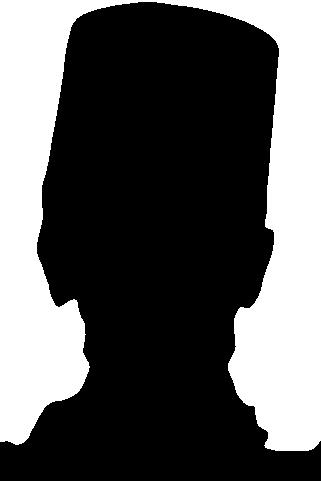}
 \\[1mm]
\multicolumn{1}{c}{\scriptsize \textit{man}} & & &\\
\hline\\
\includegraphics[width=.18\textwidth]{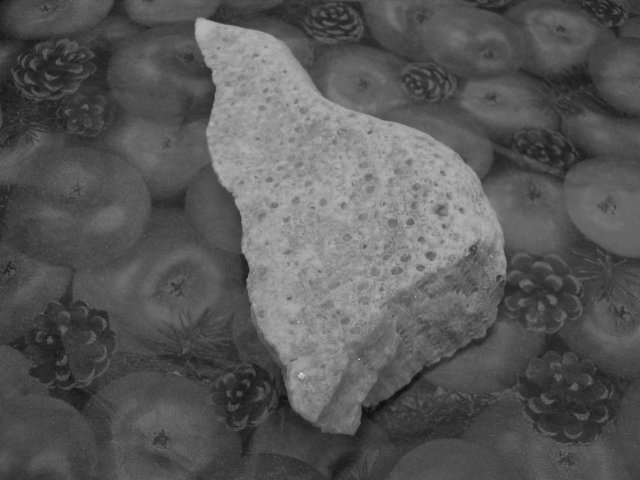} &
\includegraphics[width=.18\textwidth]{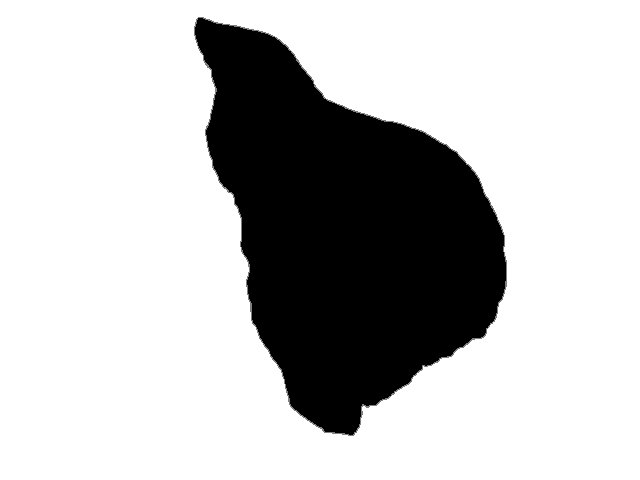}  & \includegraphics[width=.18\textwidth]{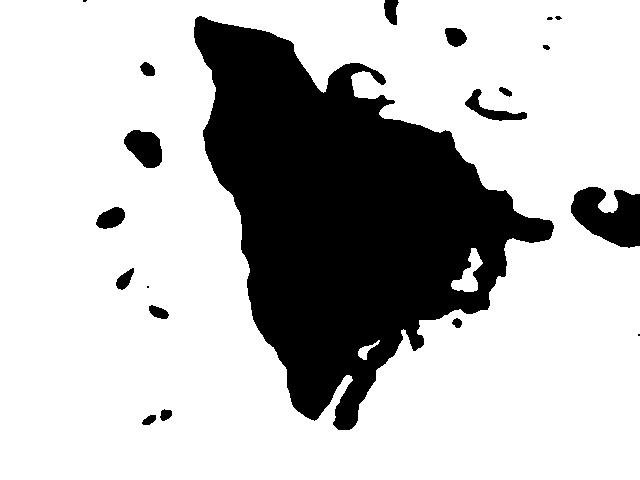} &
\includegraphics[width=.18\textwidth]{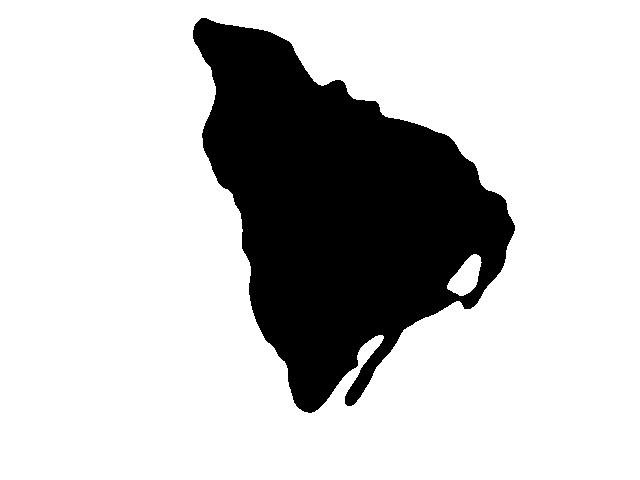}
 \\[1mm]
\multicolumn{1}{c}{\scriptsize \textit{stone}} & & &\\
\end{tabular}

\vskip 1pt
\caption{Segmentations of images with ground truth by CEN and C-TETRIS.\label{fig:gt}}
\end{center}
\end{figure*}

 First of all, in order to assess the accuracy of the C-TETRIS segmentation model, a comparison with ground truth data is presented in Fig. \ref{fig:gt}. The quality of the produced segmentations confirms the greater ability of C-TETRIS with respect to CEN in separating foreground objects from the background, especially on the \textit{flowerbed} and \textit{stone} images, where textured background is present. Furthermore, quantitative analysis measuring the similarity between the segmented images and the corresponding ground truth is given in Table \ref{tab:segerror}. The segmentation errors have been evaluated using four traditional measures\footnote{The software used for the four measures of segmentation error is available at: \url{https://people.eecs.berkeley.edu/~yang/software/lossy segmentation/.}}. The Rand Index (RI) \cite{bib:RI} counts the fraction of pairs of pixels whose labellings are consistent between the computed segmentation and the ground truth, the Global Consistency Error (GCE) \cite{bib:GCE} measures the distance between two segmentations assuming that one segmentation must be a refinement of the other, the Variation of Information (VI) \cite{bib:VI} computes the distance between two segmentations as the average conditional entropy of one segmentation given the other, and the Boundary Displacement Error (BDE) \cite{bib:BDE} computes the average boundary pixels displacement error between two segmented images\footnote{The error of one boundary pixel is defined as its distance from the closest pixel in the other boundary image.}. As we can note in Table \ref{tab:segerror}, the segmentations produced by C-TETRIS,  have smaller values of CGE, VI, and BDE, than the ones produced by CEN, as well as they present the highest values of the RI measures, showing a greater consistency with the corresponding ground truth in the partitioning of foreground objects from the background.
\begin{table}[!h]
    \centering
    \begin{tabular}{|c|c|c|c|c|c|}
    \hline 
     image & model & RI & GCE & VI & BDE \\[1mm]
     \hline 
 \multirow{2}{*}{\textit{flowerbed}}& CEN & 9.5843e-01 & 3.8905e-02 & 2.6408e-01 & 5.5895e+01 \\
 & C-TETRIS & 9.7375e-01 & 2.5570e-02 & 1.8955e-01 & 4.5199e+00 \\[1mm]
  \hline
 \multirow{2}{*}{\textit{man}}& CEN & 7.8042e-01 & 2.1254e-01 & 1.0630e+00 & 1.9497e+01\\
 & C-TETRIS   & 7.9876e-01 & 1.7726e-01 & 8.8154e-01 & 1.1459e+01\\[1mm]
\hline
  \multirow{2}{*}{\textit{stone}}& CEN & 8.9000e-01 & 1.0601e-01 & 6.1420e-01 & 2.4565e+01 \\
 & C-TETRIS & 9.1917e-01 & 7.4804e-02 & 4.5966e-01 & 1.1604e+01\\[1mm]
\hline 
    \end{tabular}
    \caption{Measures of segmentation error produced by CEN and C-TETRIS on figures displayed in Fig. \ref{fig:gt}.}
    \label{tab:segerror}
\end{table}
%
\subsection{Results on smooth images\label{subsec:smoothimages}}
\begin{figure*}[ht!]
\medskip
\begin{center}
\newcolumntype{C}{>{\centering\arraybackslash} m{.18\textwidth} }
\begin{tabular}{CCCCC}
\small original image  &
\small cartoon  &
\small CEN &
\small CEN on cartoon&
\small C-TETRIS \\
\hline\\
\includegraphics[width=.18\columnwidth]{airplane2.jpg} &
\includegraphics[width=.18\columnwidth]{airplane2U_SBit0.jpg} &
\includegraphics[width=.18\columnwidth]{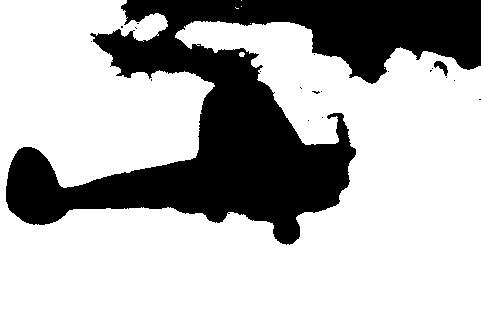} &
\includegraphics[width=.18\columnwidth]{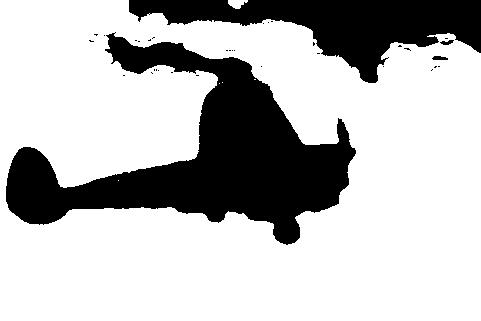} &
\includegraphics[width=.18\columnwidth]{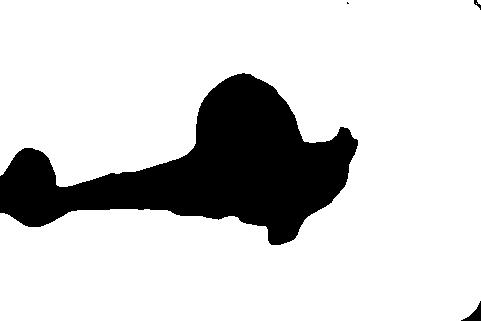} \\[1mm]
\multicolumn{1}{c}{\scriptsize \textit{airplane}} & &  \\
\hline \\%
\includegraphics[width=.18\columnwidth]{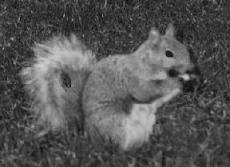} &
\includegraphics[width=.18\columnwidth]{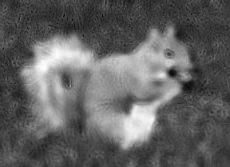} &
\includegraphics[width=.18\columnwidth]{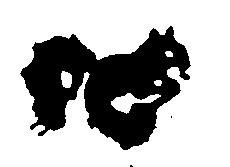} &
\includegraphics[width=.18\columnwidth]{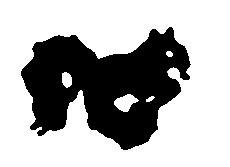} &
\includegraphics[width=.18\columnwidth]{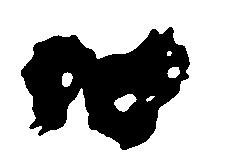} \\[1mm]
\multicolumn{1}{c}{\scriptsize \textit{squirrel}} & &  \\[1mm]
\hline \\
\includegraphics[width=.18\columnwidth]{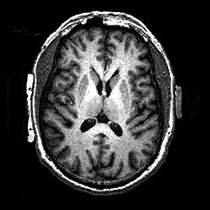} &
\includegraphics[width=.18\columnwidth]{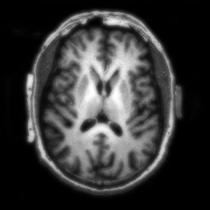} &
\includegraphics[width=.18\columnwidth]{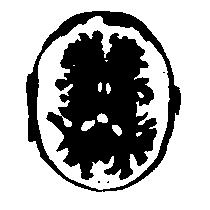} &
\includegraphics[width=.18\columnwidth]{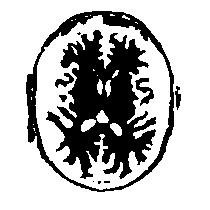} &
\includegraphics[width=.18\columnwidth]{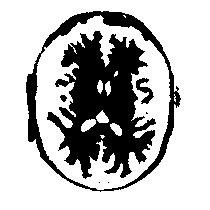} \\[1mm]
\multicolumn{1}{c}{\scriptsize \textit{brain}} & &  \\[1mm]
\hline \\
\includegraphics[width=.18\columnwidth]{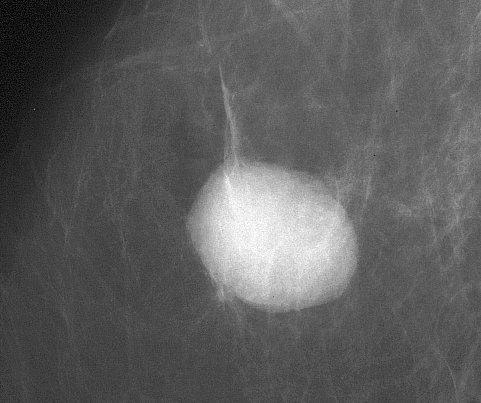} &
\includegraphics[width=.18\columnwidth]{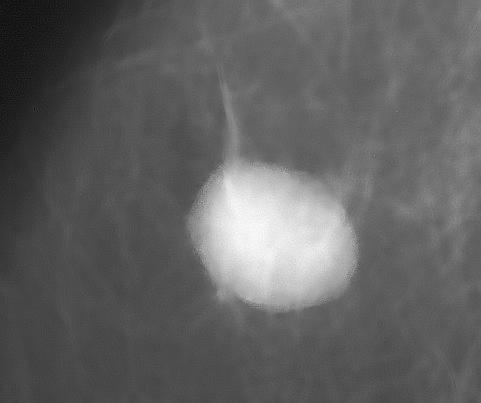} &
\includegraphics[width=.18\columnwidth]{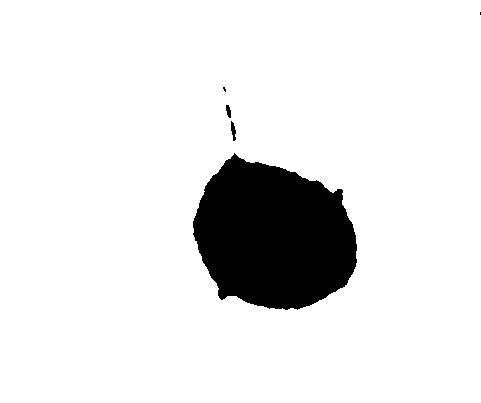} &
\includegraphics[width=.18\columnwidth]{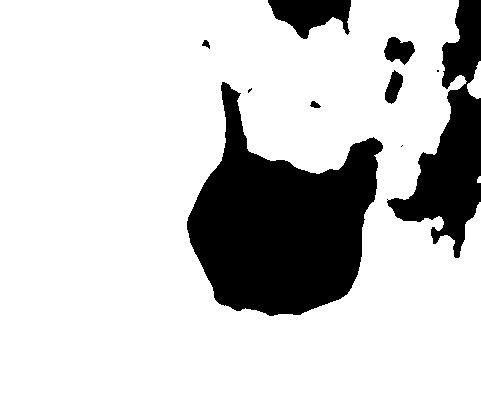} &
\includegraphics[width=.18\columnwidth]{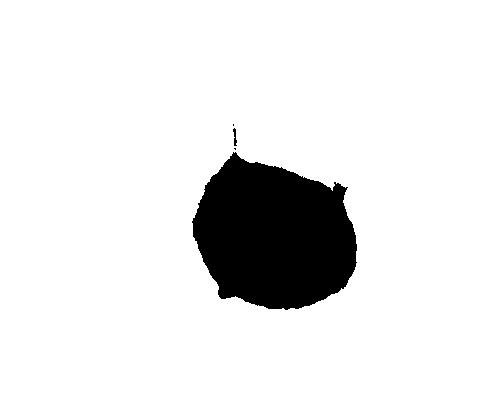} \\[1mm]

\multicolumn{1}{c}{\scriptsize \textit{ultrasound}} & &  \\[1mm]
\end{tabular}
\caption{Segmentations of smooth images by CEN and C-TETRIS. The results of the segmentation produced by CEN on the cartoon part of the image are also shown. \label{fig:CENvsCTETRIS}}
\end{center}
\end{figure*}

In Fig. \ref{fig:CENvsCTETRIS}, we show a comparison between C-TETRIS and CEN on the segmentation of the set of smooth images.
For the sake of completeness we report also the segmentation results produced by CEN on the cartoon of the images. In general, the segmentations produced by C-TETRIS are comparable with or better than the ones produced by CEN. The segmentation of \textit{airplane} shows the great effectiveness of the proposed model to separate accurately a non-uniform background from the object, due to the ability of C-TETRIS to remove the remaining texture in the cartoon, as showed in Fig. \ref{fig:CTevolution}. We note that in general there are no significant differences in the quality of the segmentation results between CEN applied to the original image and CEN applied to the cartoon. However, in the case of \textit{ultrasound}
the segmentation on the cartoon produces unreliable result,  due to the loss of contrast introduced by decomposition. In Table \ref{tab:contrast}, two global metrics are listed to measure the contrast between the given image and its cartoon. In particular we used $$m_1 = f_{max} - f_{min}$$ and the Michelson formula \cite{bib:Michelson1927}: $$m_2= (f_{max} - f_{mean})/(f_{max} + f_{mean})$$ 
where $f_{max}$, $f_{min}$ and $f_{mean}$ are the maximum, the minimum and the mean value respectively of the given image intensity. We can note that the cartoon part of \textit{ultrasound} shows the largest reduction of the both metrics with respect to the original image. 
\begin{table}[ht!]
    \centering
    \begin{tabular}{c|c|c}
        \multirow{2}{*}{image} & \multicolumn{2}{c}{contrast metrics} \\ 
        & $m_1$ & $m_2$ \\ \hline
        \textit{airplane} &  $1.$ & $0.257$ \\
        cartoon & $0.996$ & $ 0.264$\\ \hline \hline
        \textit{squirrel} & $1.$ & $0.423$ \\
        cartoon & $0.992$ & $0.442$ \\\hline \hline
        \textit{brain} & $1.$ &  $0.649$ \\
        cartoon & $1.$ & $0.645$ \\\hline \hline
        \textit{ultrasound} & $1.$ & $0.465$ \\
        cartoon & $0.953$ & $0.445$\\
                     
    \end{tabular}
    \caption{Global metrics of the image contrast (defined in section \ref{subsec:smoothimages}) evaluated on the set of smooth images and their cartoon part displayed in Fig. \ref{fig:CENvsCTETRIS}.}
    \label{tab:contrast}
\end{table}

\subsection{Results on noisy images}\label{subsec:noiseimages}
\begin{figure*}[t]
\medskip
\begin{center}
\newcolumntype{C}{>{\centering\arraybackslash} m{.24\textwidth} }
\begin{tabular}{CCC}
\multicolumn{1}{C}{\small original image }  &
\multicolumn{1}{C}{\small CEN }  &
\multicolumn{1}{C}{\small C-TETRIS }\\
\hline\\[-3mm]
\includegraphics[width=.18\textwidth]{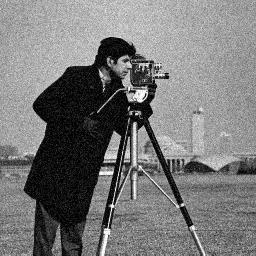} &
\includegraphics[width=.18\textwidth]{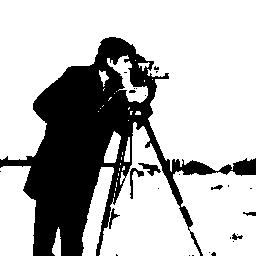} &
\includegraphics[width=.18\textwidth]{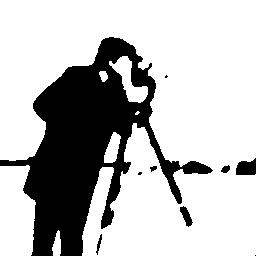} \\
\multicolumn{1}{C}{\scriptsize \texttt{gaussian noise (SNR = 20)}} & & \\
\hline\\[-3mm]
\includegraphics[width=.18\textwidth]{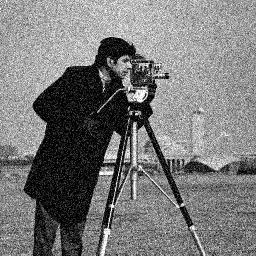} &
\includegraphics[width=.18\textwidth]{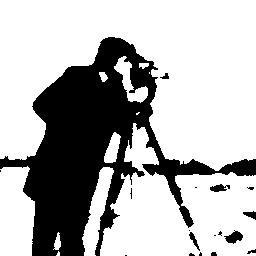} &
\includegraphics[width=.18\textwidth]{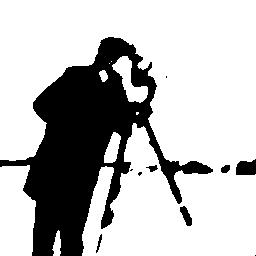} \\
\multicolumn{1}{C}{\scriptsize \texttt{gaussian noise (SNR = 15)}} & &  \\
\hline\\[-3mm]
\includegraphics[width=.18\textwidth]{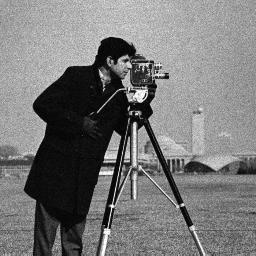} &
\includegraphics[width=.18\textwidth]{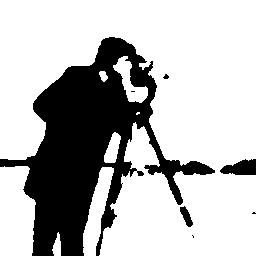} &
\includegraphics[width=.18\textwidth]{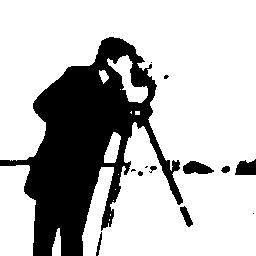} \\
\multicolumn{1}{C}{\scriptsize \texttt{poissonian noise (SNR = 35)}} & &  \\
\hline\\[-3mm]
\includegraphics[width=.18\textwidth]{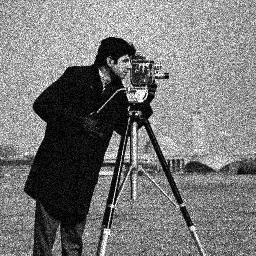} &
\includegraphics[width=.18\textwidth]{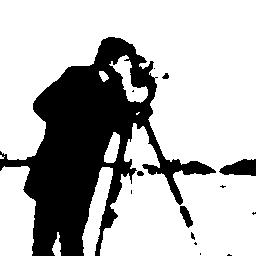} &
\includegraphics[width=.18\textwidth]{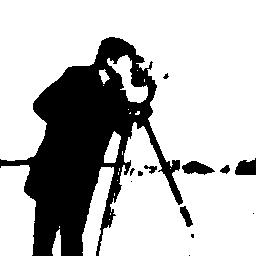} \\
\multicolumn{1}{C}{\scriptsize \texttt{poissonian noise (SNR = 30) }} & &  \\
\hline\\[-3mm]
\includegraphics[width=.18\textwidth]{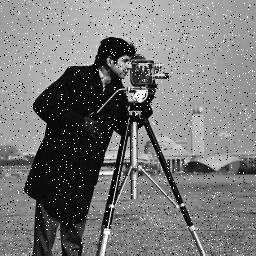} &
\includegraphics[width=.18\textwidth]{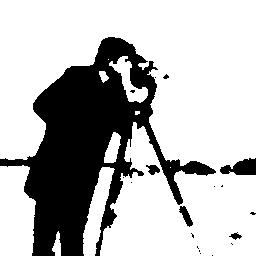} &
\includegraphics[width=.18\textwidth]{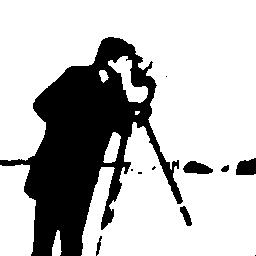} \\
\multicolumn{1}{C}{\scriptsize \texttt{salt \& pepper (5$\%$)}} & &  \\
\hline\\[-3mm]
\includegraphics[width=.18\textwidth]{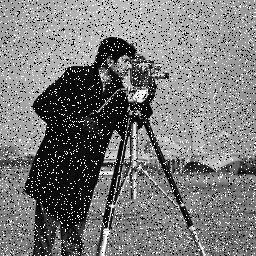} &
\includegraphics[width=.18\textwidth]{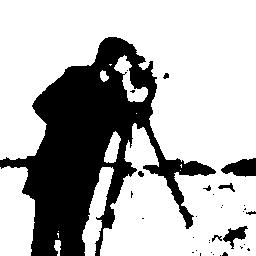} &
\includegraphics[width=.18\textwidth]{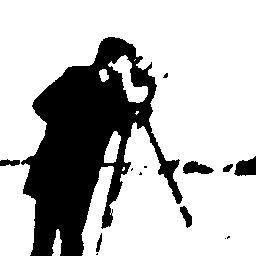} \\
\multicolumn{1}{C}{\scriptsize \texttt{salt \& pepper  (15$\%$)}} & &  \\
\end{tabular}

\vskip 1pt
\caption{Segmentations of \textit{cameraman} with different sources and noise levels by CEN and C-TETRIS. Gaussian and Poissonian noise are applied with different SNR values, whereas salt and pepper noise is added with different percentages (see section \ref{subsec:noiseimages} for details). 
\label{fig:noise}}
\end{center}
\end{figure*}

In Fig. \ref{fig:noise} a comparison between C-TETRIS and CEN on the set of noisy images is shown. The \textit{cameraman} image was corrupted by different source of noise using the MATLAB \texttt{imnoise} function. In detail: the option \texttt{`gaussian'} was used with different values for the standard deviation to obtain images affected by Gaussian noise with signal-to-noise ratio (SNR) equal to 20 and 15, respectively; by rescaling the pixels of the original image and using the option \texttt{`poisson'} we obtained images affected by Poisson noise with SNR equal to 35 and 30, respectively; finally, the option \texttt{`salt \& pepper'} was used to create images affected by impulsive noise on 5\% and 15\% of the pixels.
We note that C-TETRIS is more accurate in separating background and foreground, especially when the noise level increases. In this case, indeed, the noise is recognised as texture part and classified as foreground.

\subsection{Results on textural images\label{subsec:texturalimages}}
\begin{figure*}[!ht]
\medskip
\begin{center}
\newcolumntype{C}{>{\centering\arraybackslash} m{.20\textwidth} }
\begin{tabular}{CCCC}
\multicolumn{1}{c}{\small original image }  &
\multicolumn{1}{c}{\small C-TETRIS }  &
\multicolumn{1}{c}{\small SpAReg } &
\multicolumn{1}{c}{\small HTB }\\
\hline\\
\includegraphics[width=.18\textwidth]{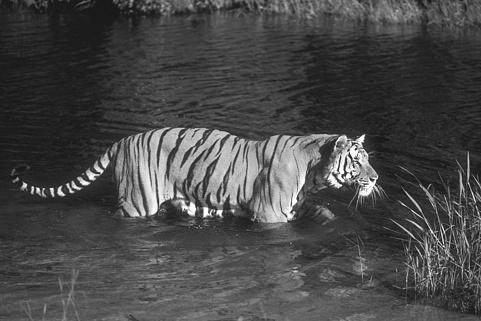} &
\includegraphics[width=.18\textwidth]{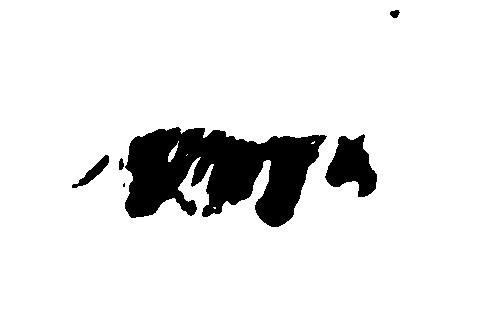} & \includegraphics[width=.18\textwidth]{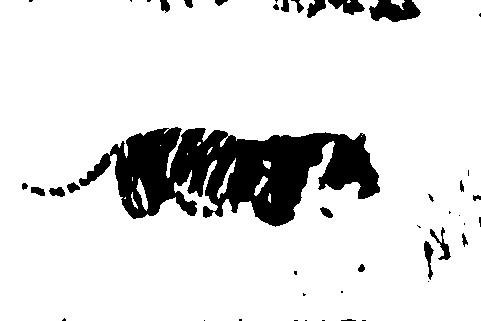} & \includegraphics[width=.18\textwidth]{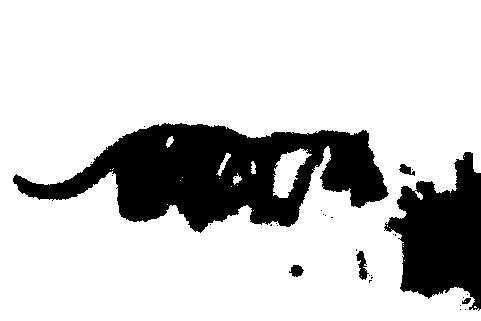}
 \\[1mm]
\multicolumn{1}{c}{\scriptsize \textit{tiger}} & & &\\
\hline\\
\includegraphics[width=.18\textwidth]{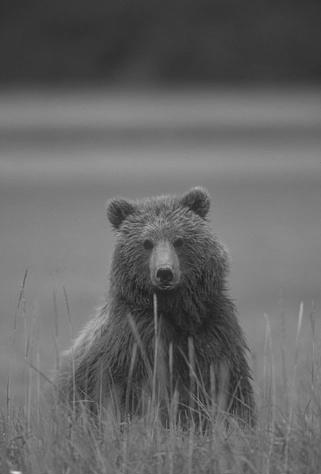} &
\includegraphics[width=.18\textwidth]{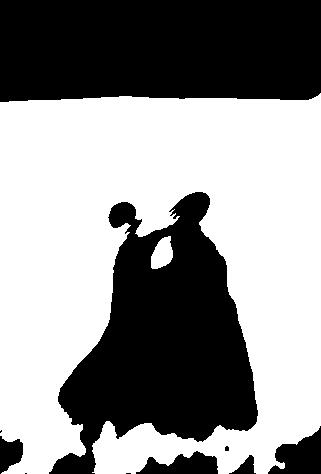} &
\includegraphics[width=.18\textwidth]{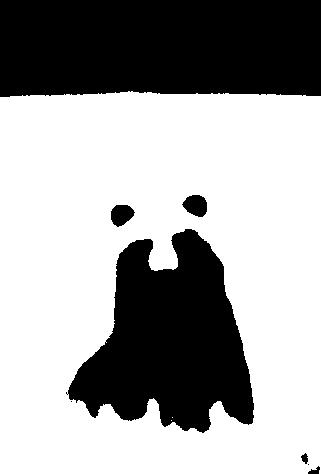} &
\includegraphics[width=.18\textwidth]{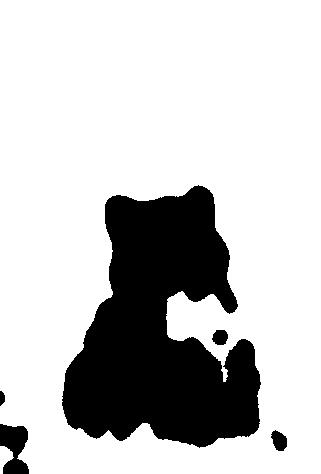} \\[1mm]
\multicolumn{1}{c}{\scriptsize \textit{bear}} & & &\\[1mm] 
\hline\\
\includegraphics[width=.18\textwidth]{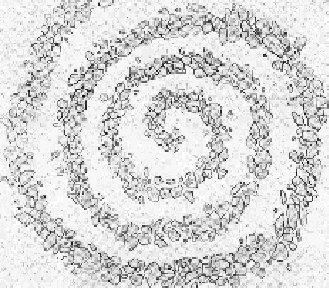} &
\includegraphics[width=.18\textwidth]{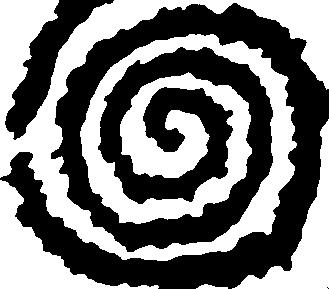} &
\includegraphics[width=.18\textwidth]{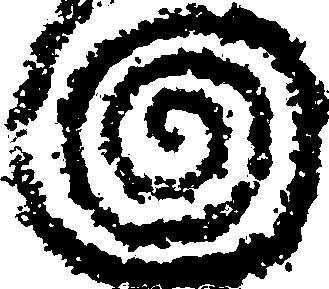} &
\includegraphics[width=.18\textwidth]{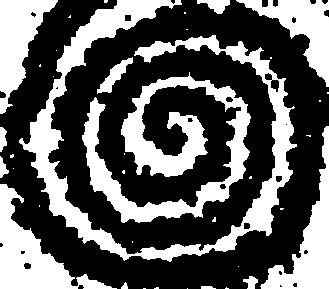} \\[1mm]
\multicolumn{1}{c}{\scriptsize \textit{spiral}} & & &\\[1mm] 
\end{tabular}

\vskip 1pt
\caption{Segmentations of textural images by C-TETRIS, SpAReg, and HTB. \label{fig:texture}}
\end{center}
\end{figure*}

Here we analyze the results of the C-TETRIS model on images containing textural components which require a two-region segmentation. We compared C-TETRIS with the Spatially Adaptive Regularization (SpAReg) model \cite{bib:Antonelli2020Adaptive}, which modifies the CEN model as follows:

\begin{equation} \label{eq:spareg}
\begin{array}{rl}
    \underset{f}{\min} &  \displaystyle \sum_{i,j} \left(\vert \nabla_x f \vert_{i,j} + \vert \nabla_y f \vert_{i,j} + \lambda_{i,j} \, (r^\top  f)_{i,j} \right) \\
   \stt             & 0 \le f \le 1
\end{array}
\end{equation}
where each entry of the matrix $\Lambda=(\lambda_{i,j})$ weighs the pixel $(i,j)$ according to local texture information as follows:
\begin{equation}
\label{eq:ctdlambda}
\lambda_{i,j} = \max \left \{ \frac{\lambda_{min}}{\lambda_{max}}, \, 1-(\rho_\sigma)_{i,j}  \right \} \lambda_{max} \, .
\end{equation}
\noindent $(\rho_\sigma)_{i,j}$ was defined applying the equation \eqref{eq:relredrate} to the given image $\bar f$, and \mbox{$0 < \lambda_{min} < \lambda_{max} < \infty$}
is a suitable range to drive the level of regularization,
depending on the image to be segmented.
In all the tests we set $ \lambda_{min} \leq \lambda_{CEN} < \lambda_{max} $.
We also include in the comparison a well-known segmentation model designed for textural images \cite{bib:HTBmodel}, that we denote as HTB. 
While C-TETRIS and SpAReg, being based on the original CEN model, classify foreground and background as regions with different intensities, the HTB model classifies them as regions with different textural components.
In detail, it finds a contour that maximizes the KL distance between the
probability density functions of the regions inside and outside the evolving (closed) active contour, which is aimed at separating textural objects of interest from the background. The feature used to characterize the texture is based on principal curvatures $\chi$ of the intensity image considered as a 2-D manifold
embedded in $\R^3$. In detail, the objective function of the HTB model is
$$ KL(p_{in},p_{out})= \sum_{i,j} ((p_{in})_{i,j}-(p_{out})_{i,j}) \, 
  ( \log \, (p_{in})_{i,j} - \log \, (p_{out})_{i,j}), $$
 where $p_{in},p_{out}$ are the probability distribution of the texture 
  feature $\chi$ in $\Omega_{in}$ and $\Omega_{out}$, respectively, assuming
  a Gaussian distribution. We consider the implementation of HTB model provided in \cite{bib:GoldsteinBressonOsher2010}. \\
\noindent Fig. \ref{fig:texture} compares the segmentations produced by C-TETRIS with the ones produced by SpAReg and HTB, respectively. Firstly, we note that C-TETRIS outperforms both SpAReg and HTB on \textit{tiger} and \textit{spiral}, where the textural object region was well identified and separated from the background. On the \textit{bear} test image, C-TETRIS seems to identify the main object better than SpAReg; however, it mistakenly includes in the foreground region some parts of the background below the bear. Both models are outperformed by HTB, which is the only model able to include the upper part of the image in the background region. In our opinion, the inaccurate result produced by the other two models is mainly due to the inhomogeneity of the background intensity that adverses its separation from the foreground region.

\section{Conclusion}\label{sec:conclusions}
In this paper, a new model named Cartoon-Texture Evolution for Two-Region Image Segmentation
(C-TETRIS) is proposed. C-TETRIS intends to improve the CEN model, which is specifically designed for smooth images, to produce good results on a wider set of images. Indeed, starting from a rough cartoon-texture decomposition of the image to be segmented, $ \bar{f} = \bar{u} + \bar{v} $, where $\bar{u}$ and $\bar{v}$ describe the cartoon and the texture components respectively,  C-TETRIS is able to simultaneously produce
a decomposition of $\bar{u}$   as $\bar{u}=u+v$, where $v$ is enforced to be close to $\bar{v}$ and the
best approximation among all the functions that take only two values of $u$.  This is realized by combining the CEN model on $u$ and a Kullback-Leibler divergence of $v$ from $\bar{v}$. The proposed model
leads to a non-smooth constrained optimization problem solved by means of the ADMM method, for which a convergence result is provided. Numerical experiments show that, as the ADMM advances, C-TETRIS progressively subtracts from $\bar{u}$ the remaining texture, leading to a clearer distinction between background and foreground of the image. The experiments show that the proposed model is able to produce accurate two-region segmentation, comparable with or better than the one produced by state-of-the-art segmentation models, for several images also corrupted by noise or containing textural components. Furthermore, C-TETRIS seems to be independent of the type and level of noise. Future work will deal with the extension of the proposed combination of cartoon-texture decomposition and KL divergence term to more advanced image segmentation models.

\subsubsection*{Data availability}\label{sec:data}
The authors confirm that all data generated or analysed during this study are included in this article. The repositories of image tests are also reported.
\subsubsection*{Competing interests}\label{sec:interests}
The authors have no financial or proprietary interests in any material discussed in this article.

\backmatter
\bmhead{Acknowledgments}
This work was partially supported by Istituto Nazionale di Alta Matematica - Gruppo Nazionale per il Calcolo Scientifico (INdAM-GNCS), by the Italian Ministry of University and Research under grant no. PON03PE\_00060\_5, and by the VALERE Program of the University of Campania ``L. Vanvitelli''. \\We would like to thank Simona Sada (ICAR-CNR) for her technical support.

%
%
%

\bibliography{biblioVL}
\end{document}